\documentclass[11pt]{article}

\usepackage{latexsym}
\usepackage{amsfonts}
\usepackage{graphicx}
\usepackage{amsmath}
\usepackage{amssymb}

\usepackage{amscd}

\usepackage{amscd}
\usepackage{tikz} 
\usetikzlibrary{matrix,arrows}
\usepackage{tikz} 
\usetikzlibrary{matrix,arrows}

\date{  }

\title{Preservation of Trees by semidirect Products}

\author{Gabriel Zapata}

\newtheorem{theorem}{Theorem}
\newtheorem{proposition}{Proposition}

\newtheorem{lemma}{Lemma}

\def\classification{\@ifnextchar [{\@xfootnotenext}%
{\begingroup\let\protect\noexpand \xdef\@thefnmark{}\endgroup
\@footnotetext}}

\newcommand{\cat}[1]{
	\ensuremath{ 
		\mathsf{#1}
	}
}

\newcommand{\f}[1]{
	\ensuremath{ \mathrm{#1} }
}

\newcommand{\slt}[1]{
	\ensuremath{ 
		\mathrm{SL}_{_2}( #1) 
	}
}

\newcommand{\glt}[1]{
	\ensuremath{ 
		\mathrm{GL}_{_2}( #1) 
	}
}

\newcommand{\gln}[2]{
	\ensuremath{ 
		\mathrm{GL}_{_{#1} }( #2) 
	}
}

\newcommand{ \kernel }[1]{
	\ensuremath{ 
		\mathrm{Ker}\, #1
	}
}

\newcommand{\im}[1]{
	\ensuremath{ 
		\mathrm{Im}\,#1 
	}
}

\newcommand{\mar}[2]{
	\ensuremath{ 
		#1\,
			\mspace{1.5mu}
				\mathbin{
					\tikz[baseline]
						\draw[->] 
						(0pt,.5ex) -- (3.5ex,.5ex);
				}
			\mspace{2mu}
		\,#2
	}
}

\newcommand{\marto}[2]{
	\ensuremath{ 
		#1\,
			\mspace{1.5mu}
				\mathbin{
					\tikz[baseline]
						\draw[|->] 
						(0pt,.5ex) -- (3.5ex,.5ex);
				}
			\mspace{2mu}
		\,#2
	}
}

\newcommand{\idi}[1]{
	\ensuremath{ 
 			1\mspace{-4.5mu}\mathrm{l}_{_{#1}} 
	}
}

\newcommand{\mor}[3]{
	\ensuremath{ \mathrm{Hom}  _{\, #1} (#2,#3) 
	}
}

\newcommand{\con}  {
		\pmb{
				\wedge
			}	
}

\begin{document}

\maketitle

\begin{abstract}
We show that the semidirect product of a group $C$
by $A*_D B$  is isomorphic to
the free  product of $A\rtimes C$ and $B\rtimes C$
amalgamated at $D\rtimes C$, where $A$, $B$ and $C$ are arbitrary
groups.  Moreover, we apply this theorem 
to  prove that any group $G$ that acts without inversion 
on a  tree $T$ that possesses a segment $\Gamma$ for its quotient graph, 
such that, if the stabilizers of the vertex 
set $\{\,P,Q\,\}$ and  edge $y$ of a lift of  $\, \Gamma$ in
$T$ are of the form  $G_{P}\!\rtimes H$,\,
$G_{Q}\!\rtimes H$ and $G_{y}\! \rtimes H$, then $G$ is isomorphic to 
the semidirect product of $H$ by $(\,G_P \,*_{G_y} \,G_Q \,)$. 

Using our results we conclude with a non-standard verification
of the isomorphism between $\glt{\mathbb{Z}}$ and the free  product of
the dihedral  groups $D_4$ and $D_6$ amalgamated at their 
Klein-four group.
\end{abstract}

\tableofcontents

\medskip

\section{Introduction}

In this paper we provide an analysis of an interplay 
between semidirect products and free products with amalgamation
(i.e., tree products.) That is, 
we show that given any groups $A$, $B$ and $C$,
the semidirect product of a group $C$ by $A*_D B$  is isomorphic to
the free  product of the $A\rtimes C$ and $B\rtimes C$ amalgamated at
$D\rtimes C$, i.e.,  
$$
\,(\,A\underset{ _{D} }{*} B \, )\rtimes C 
	~\simeq ~ 
		(A\rtimes C) \underset{ _{D\rtimes C} }{*} ( B\rtimes C )\,.
$$
Intuitively, in the category of groups, the
semidirect product of a group  distributes (or, it is preserved) 
on the right over   
 free products with amalgamation. 
Moreover,
we show that a  group $G$ that acts without inversion on a tree
$\widetilde{\Gamma}$   such that if  $G_{P}\!\rtimes H$,  
$G_{Q}\!\rtimes H$ and $G_{y}\! \rtimes H$ are the stabilizers of the
vertex set $\{\,P,Q\,\}$ and  edge $y$ of a lift of its segment 
$\, \Gamma$ in the quotient graph $\widetilde{\Gamma}$, then
\[
G\;\simeq \;(\,G_P \,\underset{G_y}{*}\,G_Q \,)\rtimes H\,.
\]
We then give an example of the isomorphism 
$
\glt{\mathbb{Z}}
	\simeq D_4   
		*_{D_2  }\,D_6
$
using these  results.

\section{An Exact Sequence for a Tree}

\begin{lemma}
\label{lem  freeproducts semidirect sequence}
If  $A$, $B$, $C$ and $D$ are groups, then there is an exact 
sequence of form  
\begin{equation}\label{dia  exact half 2}
\begin{aligned}[c]
\begin{tikzpicture}[description/.style={fill=white}]
\matrix (m) [matrix of math nodes, row sep=0 em, column sep=2 em, 
			text height=1.5ex, text depth=0.25ex] 
	{     
		1 & A \underset{ _{D} }{*} B & 
		(A\rtimes C) \underset{ _{D\rtimes C} }{*} ( B\rtimes C )
			&  C    & 1\,.    \\ 
	}; 
\path[->,font=\scriptsize]
	(m-1-1)  edge   (m-1-2);
\path[right hook->,font=\scriptsize]
	(m-1-2)  edge node[above]   {$ \nu $}  (m-1-3);
\path[->,font=\scriptsize]
	(m-1-3)  edge node[above]   {$ \mu $}  (m-1-4)
	(m-1-4)  edge 						   (m-1-5);
\end{tikzpicture}
\end{aligned}
\end{equation}

\end{lemma}

\noindent \emph{Proof}.  
Let 
$
G = (A\underset{ _{D} }{*} B \, )
$ 
and 
$
\tilde{G}
	=(A\rtimes C) \underset{ _{D\rtimes C} }{*} ( B\rtimes C )$ 
be  pushouts of the diagrams
\begin{equation}\label{dia  free of two}
\begin{aligned}[c]
\begin{tikzpicture}[description/.style={fill=white}]
\matrix (m) [matrix of math nodes, row sep=2.5 em, column sep=2.5 em, 
			text height=1.5ex, text depth=0.25ex] 
	{ D  & A     \\
	  B  & |[gray]| G	\\ }; 
\path[->,font=\scriptsize]
	(m-1-1)  edge node[left]  {$  \iota_B $}(m-2-1)
	(m-1-1)  edge node[above] {$ \iota_A $} (m-1-2);
\path[->,gray,font=\scriptsize]
	(m-2-1)  edge node[above] {$ \beta $}   (m-2-2)
	(m-1-2)  edge node[right] {$ \alpha $}  (m-2-2);
\end{tikzpicture}
\end{aligned}
\quad\text{and}\quad
\begin{aligned}[c]
\begin{tikzpicture}[description/.style={fill=white}]
\matrix (m) [matrix of math nodes, row sep=2.5 em, column sep=1.7 em, 
			text height=1.5ex, text depth=0.25ex] 
	{ D\rtimes C &  A\rtimes C     \\
	  B\rtimes C &  |[gray]|\tilde{G}\,,	\\ }; 
\path[->,font=\scriptsize]
	(m-1-1)  edge node[left] {$ \tilde{\iota}_B $} (m-2-1)
	(m-1-1) edge node[above] {$ \tilde{\iota}_A $} (m-1-2);
\path[->,gray,font=\scriptsize]
	(m-2-1)  edge node[above] {$ \tilde{\beta} $} (m-2-2)
	(m-1-2)  edge node[right] {$ \tilde{\alpha} $} (m-2-2);
\end{tikzpicture}
\end{aligned}
\end{equation} 
where its homomorphisms are injective. In particular, $\tilde{\iota}_A$
and $\tilde{\alpha}$ are the embeddings induced from embedding $\iota_A$
and the natural embedding $\alpha$, i.e., 
\[ 
\tilde{\iota}_A \,:\, 
	\marto{d\cdot c}{\big(\,\iota_A(d)\cdot c\,\big)}
	\quad\textrm{and }\quad
		\tilde{\alpha} \,:\, \marto{a\cdot c}{(a\cdot c) \,\tilde{N} }\,, 
\]
where $\tilde{N}\trianglelefteq G$ is of the form
$
\tilde{N}:= \text{Ncl}\,
	\{\; 
		\tilde{\iota}_A(d\cdot c)\;
		\tilde{\iota}_B\big(\,(d\cdot c)^{\,-1} \,\big) 
			~|~
				d\cdot c \in D\rtimes C
	\;\}\,,
$
and \, $a\cdot c\in A\rtimes C$. The embeddings $\tilde{\beta}$ and
$\tilde{\iota}_B$ are defined in the same manner.
Now	the semidirect products $A \rtimes C$ and $B\rtimes C $ are 
described by the split extensions  
\begin{equation}\label{dia  split of two}
\begin{aligned}[c]
\begin{tikzpicture}[description/.style={fill=white}]
\matrix (m) [matrix of math nodes, row sep=0 em, column sep=3 em, 
			text height=1.5ex, text depth=0.25ex] 
	{ 
		A  & A \rtimes C & C    
	\\ }; 
\path[right hook->,font=\scriptsize]	
	(m-1-1)     edge node[above] {$ \nu_A $}       (m-1-2);
\path[->,font=\scriptsize]
	(m-1-2)  edge node[above] {$ \mu_A  $} 	    (m-1-3);
\end{tikzpicture}
\end{aligned}
\quad\text{and}\quad
\begin{aligned}[c]
\begin{tikzpicture}[description/.style={fill=white}]
\matrix (m) [matrix of math nodes, row sep=0 em, column sep=3 em, 
			text height=1.5ex, text depth=0.25ex] 
	{ 
		B  & B \rtimes C & C,   
	\\ }; 
\path[right hook->,font=\scriptsize]	
	(m-1-1)     edge node[above] {$ \nu_B $}       (m-1-2);
\path[->,font=\scriptsize]
	(m-1-2)  edge node[above] {$ \mu_B   $} 	    (m-1-3);
\end{tikzpicture}
\end{aligned}
\end{equation}
where $ \nu_A $ is the natural embedding of $A$ into $A\rtimes C$,
$\mu_A$ is the projective homomorphism defined by 
\[
\mu_A: \mar{A\rtimes C} {(A\rtimes C) \,/ A} \;\simeq\; C
	\quad\text{such~that}\quad
\mu_A: \marto{a\cdot c} {\marto{Ac}}{c}\, 
\]
The homomorphisms 
$ \nu_B $ and $\mu_B$ and  defined in the similar manner. 
Then, putting the data of  diagrams  (\ref{dia  free of two}) 
and (\ref{dia  split of two}), we get the commutative diagram
\begin{equation*} 
\begin{aligned}[c]
\begin{tikzpicture}[description/.style={fill=white}]
\matrix (m) [matrix of math nodes, row sep=1.6 em, column sep=3.3 em, 
			text height=1.5ex, text depth=0.25ex] 
	{        & A             & A\rtimes C   \\
	   D	 & |[gray]|G    & \tilde{G} 	 \\
	         &  B            & B\rtimes C     \\ 
	}; 
\path[->,font=\scriptsize]
	(m-2-1)  edge node[above] {$  \iota_A $} (m-1-2)
	(m-2-1)  edge node[below] {$ \iota_B $} (m-3-2)
	(m-1-3)  edge node[right] {$ \tilde{\alpha}$} 	(m-2-3)
	(m-3-3)  edge node[right] {$ \tilde{\beta}$} 	(m-2-3);
\path[right hook->,font=\scriptsize]
	(m-1-2)  edge  node[above] {$ \nu_A $} (m-1-3) 
	(m-3-2)  edge  node[above] {$ \nu_B $} (m-3-3);
\path[->,gray,font=\scriptsize]
	(m-1-2)  edge node[right]  {$ \alpha $} (m-2-2)
	(m-3-2)  edge node[right]  {$ \beta $}  (m-2-2);
\end{tikzpicture}
\end{aligned}
 \text{, \,where}
\begin{aligned}[c]
\begin{tikzpicture}[description/.style={fill=white}]
\matrix (m) [matrix of math nodes, row sep=1.6 em, column sep=1.1 em, 
			text height=1.5ex, text depth=0.25ex] 
	{        & \iota_A(d)             & \iota_A(d) \cdot 1_{_C}   \\
	   d	 & |[gray]|\iota_A(d) = \iota_B(d)     & 
	   			\iota_A(d) \cdot 1_{_C}= \iota_B(d)	 \cdot 1_{_C}  \\
	         &  \iota_B(d)             & \iota_B(d) \cdot 1_{_C}   \\ 
	}; 
\path[|->,font=\scriptsize]
	(m-2-1)  edge node[above] {$  \iota _A $} (m-1-2)
	(m-2-1)  edge node[below] {$  \iota_B $} (m-3-2)
	(m-1-2)  edge node[above] {$ \nu_A $}   (m-1-3)
	(m-3-2)  edge node[above] {$ \nu_B $}    (m-3-3)
	(m-1-3)  edge node[right] {$ \tilde{\alpha} $} 	(m-2-3)
	(m-3-3)  edge node[right] {$ \tilde{\beta} $} 	(m-2-3);
\path[|->,gray,font=\scriptsize]
	(m-1-2)  edge node[right]  {$  \alpha  $} 	(m-2-2)
	(m-3-2)  edge node[right]  {$  \beta  $} 	(m-2-2);
\end{tikzpicture}
\end{aligned}
\end{equation*}
since 
$\,\iota_A(d)\cdot 1_{_C}= \iota_B(d)\cdot 1_{_C}\,
	\;\Longleftrightarrow \;
\tilde{\iota}_A(d\cdot 1_{_C})=\tilde{\iota}_B(d\cdot 1_{_C})\,, 
$
i.e., the diagram commutes. In addition, since $G$ is a pushout through
$D$ by $ \iota_A$ and  $\iota_B$, there is a unique
$\nu:\mar{G}{\tilde{G}}$ such that the diagram 
\begin{equation}\label{dia  free uniqueness 1}
\begin{aligned}[c]
\begin{tikzpicture}[description/.style={fill=white}]
\matrix (m) [matrix of math nodes, row sep=1.5 em, column sep=3.5 em, 
			text height=1.5ex, text depth=0.25ex] 
	{        & |[gray]|A     &             \\
	  |[gray]|D 	 & G     & \tilde{G}	\\
	         & |[gray]|B     &               \\ 
	}; 
\path[->,gray,font=\scriptsize]
	(m-2-1)  edge node[above]  {$  \iota _A $} (m-1-2)
	(m-2-1)  edge node[below]  {$  \iota _B $} (m-3-2)
	(m-1-2)  edge[bend left=20] node[above right] 	  
				{$  \tilde{\alpha} \circ\nu_A $}   (m-2-3)
	(m-3-2)  edge[bend right=20] node[below right]   
				{$ \tilde{\beta} \circ\nu_B  $}  (m-2-3)
	(m-1-2)  edge node[right]  {$  \alpha $} (m-2-2)
	(m-3-2)  edge node[right]  {$  \beta  $} (m-2-2);
\path[->,densely dotted,font=\scriptsize]
	(m-2-2)  edge node[description]  {$ \nu $} (m-2-3);
\end{tikzpicture}
\end{aligned}
\end{equation}
commutes. In particular, let $g$ be a word in $G$. Then $g$ has 
a unique normal form
$$
g
	\,=\, 
		a_1 \,
		b_1
		\,\cdots\,
		a_n\,
		b_n N\;
$$
where each $a_i \in A$,\;
$b_i\in B$ and 
$
N:= \text{Ncl}\,
	\{\, 
		\iota_A(d)\,\iota_B(d^{\,-1}) 
			~|~
				d\in D
	\,\}\,.
$
Then defining
\begin{alignat} {2}
\nu(g)
	& \,:=\, 
		a_1 ^{\,\tilde{\alpha} \circ \nu_A}\,
		b_1 ^{\,\tilde{\beta} \circ \nu_B}
		\,\cdots\,
		a_n ^{\,\tilde{\alpha} \circ \nu_A}\,
		b_n ^{\,\tilde{\beta} \circ \nu_B} 
		\tilde{N}\;
			\,=\, 
				(a_1 \cdot 1_C )
				(b_1 \cdot 1_C )
				\,\cdots\,
				(a_n \cdot 1_C )
				(b_n \cdot 1_C )
				\tilde{N}
			& &
			\notag\\
	& \,\phantom{;}=\, 
		 a_1
		b_1
		\,\cdots\,
		a_n 
		b_n
		\tilde{N}\,,
			& &
			\notag
\end{alignat}
give us a well-defined embedding, by definition and the
uniqueness of  normal forms for free products with amalgamation.
Moreover, given  $a \in A$, we can naturally check that 
\[
\nu\circ \alpha  \,(a) 
	\,:=\,
		a\tilde{N}
	\, = \,
		\tilde{\alpha} \circ \nu_A \,(a\cdot c)
\quad\text{and}\quad
	\nu\circ \beta \,(b) 
	\,:=\,
		b\tilde{N}
	\, = \,
		\tilde{\alpha} \circ \nu_A \,(b\cdot c)\,,
\]
i.e., $\nu\circ \alpha = \tilde{\alpha}\circ \nu_A$\,
and $\,\nu\circ \ \beta = \tilde{\beta}\circ \nu_B$\,. Therefore 
 $\nu$ is a monomorphism its image is of the form
\[
\im {\nu} 
	\, = \,
		\{\, 
			a_1 \,b_1 \; \cdots a_n \, b_n \,N
			~|~
 				a_i\in A ~\con ~b_i\in B
 		\,\}				
\, = \,
		G\,.				
\]
Now we would like to extend this sequence to an exact sequence. To do
this we use the diagrams  (\ref{dia  free of two}) 
and (\ref{dia  split of two}) again to get the commutative diagram
\begin{equation*}
\begin{aligned}[c]
\begin{tikzpicture}[description/.style={fill=white}]
\matrix (m) [matrix of math nodes, row sep=1.6 em, column sep=1.3 em, 
			text height=1.5ex, text depth=0.25ex] 
	{        & A\rtimes C             & A\rtimes C/A     \\
	   D\rtimes C	 & |[gray]| \tilde{G}   &C	 \\
	         &  B\rtimes C            & B\rtimes C/B       \\ 
	}; 
\path[->,font=\scriptsize]
	(m-2-1.70)  edge node[above] {$ \tilde{\iota}_A $} (m-1-2)
	(m-2-1.290)  edge node[below] {$ \tilde{\iota}_B $} (m-3-2)
	(m-1-3)  edge node[right] {$ \wr$} 	(m-2-3)
	(m-3-3)  edge node[right] {$ \wr$} 	(m-2-3);
\path[->,font=\scriptsize]
	(m-1-2)  edge  node[above] {$ \mu_A $} (m-1-3) 
	(m-3-2)  edge  node[above] {$ \mu_B $} (m-3-3);
\path[->,gray,font=\scriptsize]
	(m-1-2)  edge node[right]  {$ \alpha $} (m-2-2)
	(m-3-2)  edge node[right]  {$ \beta $}  (m-2-2);
\end{tikzpicture}
\end{aligned}
 \text{, \,where}
\begin{aligned}[c]
\begin{tikzpicture}[description/.style={fill=white}]
\matrix (m) [matrix of math nodes, row sep=1.6 em, column sep=.7 em, 
			text height=1.5ex, text depth=0.25ex] 
	{        & \iota_A(d) \cdot c          & A \, c    \\
	   d\cdot c & |[gray]|\iota_A(d) \cdot c = \iota_B(d) \cdot c     & 
	   			                                     c = c  \\
	         &  \iota_B(d)\cdot c                 & B \, c   \\ 
	}; 
\path[|->,font=\scriptsize]
	(m-2-1)  edge node[above] {$ \tilde{\iota}_A $} (m-1-2)
	(m-2-1)  edge node[below] {$ \tilde{\iota}_B $} (m-3-2)
	(m-1-2)  edge node[above] {$ \mu_A $}   (m-1-3)
	(m-3-2)  edge node[above] {$ \mu_B $}    (m-3-3)
	(m-1-3)  edge node[right] {$ \wr $} 	(m-2-3)
	(m-3-3)  edge node[right] {$ \wr $} 	(m-2-3);
\path[|->,gray,font=\scriptsize]
	(m-1-2)  edge node[right]  {$  \alpha  $} 	(m-2-2)
	(m-3-2)  edge node[right]  {$  \beta  $} 	(m-2-2);
\end{tikzpicture}
\end{aligned}
\end{equation*}
Also, since $\tilde{G}$ is a pushout through
$D\rtimes C$ by $\tilde{\iota}_A$ and  $\tilde{\iota}_B$, then there 
exists  
a unique $\mu:\mar{\tilde{G}} {C}$ such that the diagram 
\begin{equation}\label{dia  free uniqueness 2}
\begin{aligned}[c]
\begin{tikzpicture}[description/.style={fill=white}]
\matrix (m) [matrix of math nodes, row sep=1.5 em, column sep=2.5 em, 
			text height=1.5ex, text depth=0.25ex] 
	{        & |[gray]|A\rtimes C     &                \\
	  |[gray]|D \rtimes C   	 & \tilde{G}    & C 	\\
	         & |[gray]|B \rtimes C       &               \\ 
	}; 
\path[->,gray,font=\scriptsize]
	(m-2-1)  edge node[above]  {$ \tilde{\iota}_A $} (m-1-2)
	(m-2-1)  edge node[below]  {$ \tilde{\iota}_B $} (m-3-2)
	(m-1-2)  edge[bend left=20] node[above right] 	  
    									{$  \mu_A $} (m-2-3)
	(m-3-2)  edge[bend right=20] node[below right]   
									   	{$ \mu_B  $} (m-2-3)
	(m-1-2)  edge node[right] 		  {$  \alpha  $} (m-2-2)
	(m-3-2)  edge node[right] 		   {$  \beta  $} (m-2-2);
\path[->,densely dotted,font=\scriptsize]
	(m-2-2)  edge node[description] 	   {$ \mu $} (m-2-3);
\end{tikzpicture}
\end{aligned}
\end{equation}
commutes. In particular, if $g$ is a word in $G$, then $\tilde{g}$ has a
unique normal form
\begin{alignat} {2}
\tilde{g}
	& := \,
		(a_1 \cdot c_1 )
		(b_1^{\,\prime} \cdot c_1^{\,\prime} )
		\cdots 
		(a_n \cdot c_n )
		(b_n^{\,\prime} \cdot c_n^{\,\prime} ) \tilde{N}
			& &
			\notag\\
	& \,=\, 
		 a_1 \, b_1  ^{\,\prime ~c_1} \;  
		 a_2 \,b_2^{\, \prime ~c_1c_1' c_2 \,}
		 \ldots\,
		 a_2 \,b_2^{ 
		 			\prime ~c_1 c_1' 
		 			\cdots 
		 			\;c_{_{n-1}}c_{_{n-1}}' 
		 		\,}
		 \cdot
		 (c_1c_1^\prime \cdots \,c_n c_n^{\,\prime}) \;\tilde{N}
			& &
			\notag\\
		& \, = \, 
		 a_1 \,
		 b_1   
		 \cdots
		 a_n \, 
		 b_n
		 \cdot c\, \tilde{N}
		 	& &
		 	\notag
\end{alignat}
where $a_i\cdot c_i \in A\rtimes C$,\;
$b'_i\cdot c'_i \in B\rtimes C$ and $b_i\in B$. Hence
\begin{alignat} {2}
\mu(g)
	& \,:=\, 
		 a_1^{\,\mu_A} \,
		 b_1^{\,\mu_B} \;  
		 a_2^{\,\mu_A} \,
		 b_2^{\,\mu_B}
		 \cdots
		 a_n ^{\,\mu_A}\, 
		 b_n^{\,\mu_B}
		 \cdot c\,
			& &
			\notag\\
	& \,\phantom{;}=\, 
		 c
			& &
			\notag
\end{alignat}
which is clearly a well-defined epimorphism (again, by definition and the
uniqueness of  normal forms for free products with amalgamation.) The
kernel of $\mu$ has the form
\[
\kernel {\mu} 
	\, = \,
		\{\, 
			a_1 \,b_1 \; \cdots a_n \, b_n \,N
			~|~
 				a_i\in A ~\con ~b_i\in B
 		\,\}				
\, = \,
		G\,.				
\]
Therefore, $\kernel {\mu} \, =\, \im {\nu}$.\, Moreover, given  
$a \in A$, we can  check that 
$\mu\circ \tilde{\alpha}  
	\, = \,
		\alpha \circ \mu_A \,
$
and 
$
	\,\mu\circ \tilde{\beta} 
	\,=\,
		\alpha \circ \mu_B \,.
$
Therefore
we get the diagram
\begin{equation*}
\begin{aligned}[c]
\begin{tikzpicture}[description/.style={fill=white}]
\matrix (m) [matrix of math nodes, row sep=2 em, column sep=2.4 em, 
			text height=1.5ex, text depth=0.25ex] 
	{     
		1 & A \underset{ _{D} }{*} B & 
		(A\rtimes C) \underset{ _{D\rtimes C} }{*} ( B\rtimes C )
			&  C    & 1\, ,   \\ 
	}; 
\path[->,font=\scriptsize]
	(m-1-1)  edge   (m-1-2);
\path[right hook->,font=\scriptsize]
	(m-1-2)  edge node[above]   {$ \nu $}  (m-1-3);
\path[->,font=\scriptsize]
	(m-1-3)  edge node[above]   {$ \mu $}  (m-1-4)
	(m-1-4)  edge 						   (m-1-5);
\end{tikzpicture}
\end{aligned}
\end{equation*}
which is an exact sequence.
	\hfill$\Box$

\

\section{A Preservation of a Tree by a Semidirect Product}

\begin{proposition}
\label{pro  semidirect functor}
Let $\cat{Grp}$ be the category of groups, let  $C$ be a group and  let
\[
\f{Grp}_{_{\rtimes C}}:
	\mar{\cat{Grp} } 
		{\cat{Grp}}
\]  
be the assignment defined as follows:

\begin{itemize}

\item   
$ 
\f{Grp}_{_{\rtimes C}}:
	 \marto {G} { G\rtimes C}
$
for any group $G$.

\item  
$\f{Grp}_{_{\rtimes C}}:
	\marto{\psi}
		{ \psi\rtimes \idi{C} }
$ 
for any $\psi \in \mor{}{G}{H}$ such that $\idi{C}$ is the identity
automorphism of $C$ and
\[
\psi\rtimes \idi{C} : \mar{G \rtimes C}{H \rtimes C}
	\quad \textrm{is defined by}\quad
\psi\rtimes \idi{C} \,: 
	\marto{\;g\cdot c} { g^{\psi}\cdot c }\,,
\]
where $g\in G$ and $c\in C$.

\end{itemize}
Then $\f{Grp}_{_{\rtimes C}}$ is a functor.

\end{proposition}

\noindent \emph{Proof}.  By definition, the map is well-defined on 
the class of groups and homomorphisms. Let $G$ and $H$ be groups, and 
let
$
\psi :\marto {G}{H}
$
be a homomorphism.  Suppose $g\in G$ and $c\in C$,  then 
$
\f{Grp}_{_{\rtimes C}} \big(\psi\big) \,(g\cdot c)
	\, : = \, 
		\psi\rtimes \idi{C} \,(g\cdot c)
	\, = \,
		g^\psi \cdot c\,.
$
Therefore
\begin{equation}\label{dia  semidirect functor obj}
\begin{aligned}[c]
\begin{tikzpicture}[description/.style={fill=white}]
\matrix (m) [matrix of math nodes, row sep=2.75em, column sep=.8em, 
			text height=1.5ex, text depth=0.25ex] 
	{ G      		 & \phantom{12}  & H        				\\
	  G\rtimes C	 & \phantom{22}  & H\rtimes C\,.	    	\\ }; 
\path[->,font=\scriptsize]
	(m-1-1) edge node[above]    {$ \psi $}      (m-1-3)
 	(m-2-1)edge  node[below]    {$ \psi\rtimes \idi{ C } $} (m-2-3);
\path[->,dotted,font=\scriptsize]	
	(m-1-1) edge 						   	 (m-2-1)
	(m-1-2) edge node[description]	{$ \f{Grp}_{_{\rtimes C}} $} 
	 (m-2-2.90)
    (m-1-3) edge   							 (m-2-3);
\end{tikzpicture}
\end{aligned}
\end{equation} 
In particular, if $\psi$ is the identity automorphims 
$\idi{C} : \mar{G}{G} $, then 
$
\f{Grp}_{_{\rtimes C}}( \idi{G}) 
$
is equal to
$
		\idi{G} \rtimes \idi{C}.
$
Now 
$
\idi{G} \rtimes \idi{C} (g\cdot c) 
	\, = \, 
		g\cdot c
	\,=\,
		\idi{G \rtimes C}  (g\cdot c)
$ 
and therefore
\[
\f{Grp}_{_{\rtimes C}}( \idi{G})  = \idi{G \rtimes C}\,.
\] 
Also, if
$
\varphi :\mar {H}{K}
$
, where $K$ is a group, then
\[
\f{Grp}_{_{\rtimes C}} \big(\varphi\circ \psi \big)\, (g \cdot c)
	\,=\,
		 g^{\psi\, \circ \varphi} \cdot c
	\,=\,	
		\f{Grp}_{_{\rtimes C}}(\varphi) \, (g^{\,\psi} \cdot c)
	\,=\,
		\f{Grp}_{_{\rtimes C}} (\varphi)
		\, \big(\,\f{Grp}_{_{\rtimes C}} (\psi )  (g \cdot c)\,\big)
\]
i.e., 
$
\f{Grp}_{_{\rtimes C}} (\varphi\circ \psi )
	= 
		\f{Grp}_{_{\rtimes C}} \, (\varphi) \circ
		\f{Grp}_{_{\rtimes C}} \, (\psi)\,.
$
Therefore $\f{Grp}_{_{\rtimes C}}$ is a functor.
	\hfill$\Box$
	
\

\begin{theorem}\label{the  preservation semi free}
\label{pro  preserve freeproducts under semidirect}
For any group  $C$  the functor\,   
$
\f{Grp}_{_{\rtimes C}}:
	\mar{\cat{Grp} } 
		{\cat{Grp}}
$  
\,preserves free products with amalgamation, i.e.,
\[
(\,A\underset{ _{D} }{*} B \, )\rtimes C 
	~\simeq ~ 
		(A\rtimes C) \underset{ _{D\rtimes C} }{*} ( B\rtimes C )
\]
given any group $A$, $B$ and $D$.
\end{theorem}

\noindent \emph{Proof}. Let 
$G = (A\underset{ _{D} }{*} B \, )\rtimes C$ and let 
$\tilde{G}=(A\rtimes C) \underset{ _{D\rtimes C} }{*} ( B\rtimes C )$.
By lemma \ref{lem  freeproducts semidirect sequence}
\begin{equation}\label{dia  exact 1}
\begin{aligned}[c]
\begin{tikzpicture}[description/.style={fill=white}]
\matrix (m) [matrix of math nodes, row sep=2 em, column sep=2.4 em, 
			text height=1.5ex, text depth=0.25ex] 
	{     
		\phantom{1} & A \underset{ _{D} }{*} B & 
		(A\rtimes C) \underset{ _{D\rtimes C} }{*} ( B\rtimes C )
			&  C   & \phantom{1}    \\ 
	}; 
\path[right hook->,font=\scriptsize]
	(m-1-2)  edge node[above]   {$ \nu $}  (m-1-3);
\path[->,font=\scriptsize]
	(m-1-3)  edge node[above]   {$ \mu $}  (m-1-4);
\end{tikzpicture}
\end{aligned}
\end{equation}
is exact; hence, it suffices to show that this sequence is also
a split extension. The colimit of the diagrams 
\begin{equation}\label{dia  free of two 2}
\begin{aligned}[c]
\begin{tikzpicture}[description/.style={fill=white}]
\matrix (m) [matrix of math nodes, row sep=2.5 em, column sep=2.5 em, 
			text height=1.5ex, text depth=0.25ex] 
	{ D  & A     \\
	  B  & |[gray]| G	\\ }; 
\path[->,font=\scriptsize]
	(m-1-1)  edge node[left]  {$  \iota_B $}(m-2-1)
	(m-1-1)  edge node[above] {$ \iota_A $} (m-1-2);
\path[->,gray,font=\scriptsize]
	(m-2-1)  edge node[above] {$ \beta $}   (m-2-2)
	(m-1-2)  edge node[right] {$ \alpha $}  (m-2-2);
\end{tikzpicture}
\end{aligned}
\quad\text{and}\quad
\begin{aligned}[c]
\begin{tikzpicture}[description/.style={fill=white}]
\matrix (m) [matrix of math nodes, row sep=2.5 em, column sep=1.7 em, 
			text height=1.5ex, text depth=0.25ex] 
	{ D\rtimes C &  A\rtimes C     \\
	  B\rtimes C &  |[gray]|\tilde{G}\,,	\\ }; 
\path[->,font=\scriptsize]
	(m-1-1)  edge node[left] {$ \tilde{\iota}_B $} (m-2-1)
	(m-1-1) edge node[above] {$ \tilde{\iota}_A $} (m-1-2);
\path[->,gray,font=\scriptsize]
	(m-2-1)  edge node[above] {$ \tilde{\beta} $} (m-2-2)
	(m-1-2)  edge node[right] {$ \tilde{\alpha} $} (m-2-2);
\end{tikzpicture}
\end{aligned}
\end{equation} 
are the tree products $G$ and $\tilde{G}$, respectively. These
groups are described by its commutative diagrams, where $\tilde{\iota}_A$
and $\tilde{\alpha}$ are the embeddings induced from embedding $\iota_A$
and the natural embedding $\alpha$ (as described in the proof of
lemma \ref{lem  freeproducts semidirect sequence}.) The embeddings
$\tilde{\beta}$ and $\tilde{\iota}_B$ are defined in the same manner.
Now	the semidirect products $A \rtimes C$ and $B\rtimes C $ are 
described by the split extensions  
\begin{equation}\label{dia  split of two 2}
\begin{aligned}[c]
\begin{tikzpicture}[description/.style={fill=white}]
\matrix (m) [matrix of math nodes, row sep=0 em, column sep=3 em, 
			text height=1.5ex, text depth=0.25ex] 
	{ 
		A  & A \rtimes C & C    
	\\ }; 
\path[right hook->,font=\scriptsize]	
	(m-1-1)     edge node[above] {$ \nu_A $}       (m-1-2);
\path[->,font=\scriptsize]
	(m-1-2.10)  edge node[above] {$ \mu_A  $} 	    (m-1-3.158);
\path[<-,densely dotted,font=\scriptsize]
	(m-1-2.355) edge node[below] {$ \tau_A  $} 	(m-1-3.192);
\end{tikzpicture}
\end{aligned}
\quad\text{and}\quad
\begin{aligned}[c]
\begin{tikzpicture}[description/.style={fill=white}]
\matrix (m) [matrix of math nodes, row sep=0 em, column sep=3 em, 
			text height=1.5ex, text depth=0.25ex] 
	{ 
		B  & B \rtimes C & C,   
	\\ }; 
\path[right hook->,font=\scriptsize]	
	(m-1-1)     edge node[above] {$ \nu_B $}       (m-1-2);
\path[->,font=\scriptsize]
	(m-1-2.10)  edge node[above] {$ \mu_B   $} 	    (m-1-3.160);
\path[<-,densely dotted,font=\scriptsize]
	(m-1-2.355) edge node[below] {$ \tau_B  $} 	(m-1-3.192);
\end{tikzpicture}
\end{aligned}
\end{equation}
where $ \nu_A $ is the natural embedding, $\mu_A$ is the projective
homomorphism and $\tau_A$ is the transversal homomorphism of the split
extension. In particular, for any $g\in A\rtimes C$ there are 
unique $a\in A$ and $c\in C$ such that $g = a\cdot c$. Then 
$
Ag
	\,=\,
		A(a\cdot c)
	\,=\,
		A\,c
$
and 
\[
\tau_A: \mar{ C\simeq(A\rtimes C)\,/A } { A\rtimes C } 
	\quad\text{is defined by}\quad
\tau_A: \marto{c\simeq A\,g} { 1_A\cdot c } \,,
\]
which satisfies $\mu_A\circ\tau_A = \idi{A}$.
The homomorphisms 
$ \nu_B $, $\mu_B$ and $\tau_B$ are defined in the similar manner, as
before. Let us define a homomorphism $\tau:\mar{C} {\tilde{G}}$ by 
the natural embedding  $\tau:\mar{ c } { c\tilde{N} }$. Then, putting 
the data of  the diagrams  (\ref{dia  free of two}) and 
(\ref{dia  split of two}),  we get the commutative diagram 
\begin{equation*}\label{dia splitting}
\begin{aligned}[c]
\begin{tikzpicture}[description/.style={fill=white}]
\matrix (m) [matrix of math nodes, row sep=1.4 em, column sep=2.8 em, 
			text height=1.5ex, text depth=0.25ex] 
	{           & |[gray]|A\rtimes C   \\
	     C      & \tilde{G}     \\
	            & |[gray]|B\rtimes  C     \\ 
	}; 
\path[->,gray,font=\scriptsize]
	(m-2-1)  edge[bend left=20] node[above] {$ \tau_A$} 	(m-1-2)
	(m-2-1)  edge[bend right=20] node[below] {$ \tau_B$} 	(m-3-2);
\path[->,gray,font=\scriptsize]
	(m-1-2)  edge node[right]  {$ \alpha $} (m-2-2)
	(m-3-2)  edge node[right]  {$ \beta $}  (m-2-2);
\path[->,densely dotted, font=\scriptsize]
	(m-2-1)  edge node[description]  {$ \tau $}  (m-2-2);
\end{tikzpicture}
\end{aligned}
\end{equation*}
In particular, given $c\in C$, then 
$
\mu\circ \tau(c)
  \,:=\, 
		\mu (c\,\tilde{N})
 	\,=\, 
		 c\,,
$
i.e.,   $ \mu\circ \tau = \idi{C}$. Therefore, $\tau$ is a 
transversal homomorphism of the extension represented by
diagram  (\ref{dia  exact 1}). Hence,
\begin{equation}\label{dia  exact 2}
\begin{aligned}[c]
\begin{tikzpicture}[description/.style={fill=white}]
\matrix (m) [matrix of math nodes, row sep=2 em, column sep=2.4 em, 
			text height=1.5ex, text depth=0.25ex] 
	{     
		1 & A \underset{ _{D} }{*} B & 
		(A\rtimes C) \underset{ _{D\rtimes C} }{*} ( B\rtimes C )
			&  C   & 1    \\ 
	}; 
\path[->,font=\scriptsize]
	(m-1-1)  edge   (m-1-2)
	(m-1-4)  edge   (m-1-5);
\path[right hook->,font=\scriptsize]
	(m-1-2)  edge node[above]   {$ \nu $}  (m-1-3);
\path[->,font=\scriptsize]
	(m-1-3)  edge node[above]   {$ \mu $}  (m-1-4);
\end{tikzpicture}
\end{aligned}
\end{equation}
is a split extension; consequently
$
\,(\,A\underset{ _{D} }{*} B \, )\rtimes C 
	~\simeq ~ 
		(A\rtimes C) \underset{ _{D\rtimes C} }{*} ( B\rtimes C )\,.
$
	\hfill$\Box$

\begin{theorem}
Let  $G$ act without inversion on a tree $\widetilde{\Gamma}$ and
let $\Gamma = G \backslash \widetilde{\Gamma}$  denote its factor graph.
If $\Gamma$ is a segment such that $G_{P}\!\rtimes H$, 
$G_{Q}\!\rtimes H$ and $G_{y}\! \rtimes H$ are  the stabilizers of the
vertex set $\{\,P,Q\,\}$ and  edge $y$ of a lift of 
$\, \Gamma$ in $\widetilde{\Gamma}$,  then
\[
G\;\simeq \;(\,G_P \,\underset{G_y}{*}\,G_Q \,)\rtimes H\,.
\]
\end{theorem}

\noindent \emph{Proof}. The canonical  homomorphism
\[
\varphi: 
	\mar{
			(\,G_P \rtimes H\,)\,
				\underset{G_y\,\rtimes H}{*}
			\,(\,G_Q\rtimes H \,) 
		}
	    { G }
\]
is an isomorphism, \emph{a priori}. Therefore 
$\,
G\;\simeq \;(\,G_P \,\underset{G_y}{*}\,G_Q \,)\rtimes H
\,$
by theorem  \ref{the  preservation semi free}.
	\hfill	$\square$

\section{A Quick Application to $\glt{\mathbb{Z}}$}

The \emph{general linear group}
$
\gln{n}{R}
$
over a commutative $R$-module $M$ consists of
the set of invertible linear operators over $M$. We are 
particularly interested  in the group
$ 
\glt{\mathbb{Z}}\,,
$  
which  has a representation 
\[
\glt{\mathbb{Z}}
	\,:=\,
		\{\, A\in \mathrm{M}_{_{2\times 2}}
			~|~ \det{A} = \pm 1 
		\,\}
\]
since any invertible matrix $A\in \mathrm{M}_{_{2\times 2}}$
is the form 
$
A^{-1} \,=\, \left(\frac{1}{\det A}\right) \mathrm{adj}\, A\,,
$
where $\mathrm{adj}\,A$ is the \emph{classical adjoint}
of $A$ and $\det A$ is the \emph{determinant} of $A$.  
The \emph{special linear group} of $ 2\times 2$
matrices, denoted by $\slt{\mathbb{Z}}$, is a subgroup  
$\glt{\mathbb{Z}}$ with representation
$
\slt{\mathbb{Z}}
	\,:=\,
		\{\, A\in \mathrm{M}_{_{2\times 2}}
			~|~ \det{A} =  1 
		\,\}\,.
$
The theory of group actions on trees tells us that
\[
\slt{\mathbb{Z}}
	\;\simeq\;
		\mathbb{Z}_4\, \underset{ \;\mathbb{Z}_2  }{*} \,\mathbb{Z}_6\,.
\]
We now use the above result to relate  $\glt{\mathbb{Z}}$ with 
$\slt{\mathbb{Z}}$ by way of theorem  
\ref{the  preservation semi free}.

\begin{theorem}\label{prop  gl2 as a amalgamation}
$\f{GL}_2(\mathbb{Z})\simeq D_4\, \underset{ _{ D_2 } }{*}\,D_6$.
\end{theorem}

\noindent\emph{Proof}. The sequence of groups
\begin{equation}\label{dia  sl2 exact sequence}
\begin{aligned}[c]
\begin{tikzpicture}[description/.style={fill=white}]
\matrix (m) [matrix of math nodes, row sep=2 em, column sep=2 em, 
			text height=1.5ex, text depth=0.25ex] 
	{     
		1 & \slt{ \mathbb{Z} }  & 
			\glt{ \mathbb{Z} }  & 
				\mathbb{Z}_2    & 1    \\ 
	}; 
\path[->,font=\scriptsize]
	(m-1-1)  edge   (m-1-2)
	(m-1-3)  edge node[above]   {$ \mathrm{det} $}  (m-1-4)
	(m-1-4)  edge 						   (m-1-5);
\path[right hook->,font=\scriptsize]
	(m-1-2)  edge     (m-1-3);
\end{tikzpicture}
\end{aligned}
\end{equation}
is exact, since $\slt{\mathbb{Z}}= \kernel (\det)$. The map 
\[
\tau:\mar{\mathbb{Z}_2 }{ \glt{\mathbb{Z}} }
\quad\text{defined~by}\quad
\tau:\marto{-1} 
{
\left( 
	 	\begin{smallmatrix}	 
		 0 & 1				\\			       
		 1 & 0 			     \\
		\end{smallmatrix}
\right)
}
\]
is  clearly a well-defined homomorphism. In particular, since the
$(\det \circ \,\tau) = \idi{A}$,
the extension in diagram
(\ref{dia  sl2 exact sequence}) is a split extension; hence 
$\glt{\mathbb{Z}}\simeq\slt{ \mathbb{Z} }\rtimes \mathbb{Z}_2$.
Thus
\begin{alignat} {2}
\slt{\mathbb{Z}}\rtimes \mathbb{Z}_2
	& \;\simeq\;
		(\,\mathbb{Z}_4 \underset{ _{\mathbb{Z}_2 } }{*} \mathbb{Z}_6\,)
		\rtimes \mathbb{Z}_2
			& &
			\quad \text{--- \emph{a priori};}
			\notag\\
& \;\simeq\;
		(\,\mathbb{Z}_4\rtimes \mathbb{Z}_2\,)
		 \underset{  \mathbb{Z}_2\rtimes \mathbb{Z}_2  }{*}
		  (\,\mathbb{Z}_6\rtimes \mathbb{Z}_2\,)
			& &
			\quad\text{--- by theorem 
				\ref{pro  preserve freeproducts under semidirect}; and}
			\notag\\
& \;\simeq\;
		  D_4 
		  \underset{ _{ D_2 } }{*}
		  D_6 \,.
			& &
			\notag
\end{alignat}
Therefore 
$\glt{\mathbb{Z}}
	\simeq D_4   
		\underset{ _{ D_2 } }{*}D_6\,.
$ 
	\hfill$\Box$

\vspace*{15pt}


\medskip




\end{document}